\documentclass[reqno]{amsart}
\usepackage{float}
\usepackage{amsmath}
\usepackage{graphicx}
\usepackage{latexsym}
\usepackage{amsfonts}
\usepackage{amssymb}
\usepackage{color}
\theoremstyle{plain}
\newtheorem{theorem}{Theorem}
\newtheorem{corollary}[theorem]{Corollary}

\theoremstyle{definition}

\newtheorem{remark}[theorem]{Remark}
\newtheorem*{remark*}{Remark}

\newcommand{\rd}{\mathbb R^d}
\newcommand{\pr}{\mathbf P}

\begin{document}
\title[Invariance principles for random walks in cones]{Invariance principles for random walks in cones}
\author[Duraj]{Jetlir Duraj}
\address{Department of Economics, Harvard University, USA }
\email{duraj@g.harvard.edu}

\author[Wachtel]{Vitali Wachtel}
\address{Institut f\"ur Mathematik, Universit\"at Augsburg, 86135 Augsburg, Germany}
\email{vitali.wachtel@mathematik.uni-augsburg.de}

\begin{abstract}
We prove invariance principles for a multidimensional random walk conditioned to stay in a cone.
Our first result concerns convergence towards the Brownian meander in the cone. Furthermore,
we prove functional convergence of $h$-transformed random walk to the corresponding $h$-transform of the 
Brownian motion. Finally, we prove an invariance principle for bridges of a random walk in a cone.
\end{abstract}


\maketitle
{\scriptsize
}

\section{Introduction, main results and discussion}
Consider a random walk $\{S(n),n\geq1\}$ on $\rd$, $d\geq1$, where
$$
S(n)=X(1)+\cdots+X(n)
$$
and $\{X(n), n\geq1\}$ is a family of independent copies of a random variable $X=(X_1,X_2,\ldots,X_d)$.
Denote by $\mathbb{S}^{d-1}$ the unit sphere of $\rd$ and $\Sigma$ an open and connected subset of
$\mathbb{S}^{d-1}$. Let $K$ be the cone generated by the rays emanating from the origin and passing
through $\Sigma$, i.e. $\Sigma=K\cap \mathbb{S}^{d-1}$.

Let $\tau_x$ be the exit time from $K$ of the random walk with starting point $x\in K$, that is,
$$
\tau_x=\inf\{n\ge 1: x+S(n)\notin K\}.
$$

Denisov and Wachtel \cite{DW15} have determined the asymptotic behaviour of $\mathbf{P}(\tau_x>n)$ and
have proven limit theorems (integral and local) for $S(n)$ on the event $\{\tau_x>n\}$.

The main purpose of this paper is to continue the study of random walks in cones and to derive invariance
principles for various functionals of random walks constrained to stay in a cone. The proofs in \cite{DW15} are based in the strong approximation of
random walks by the Brownian motion. Therefore, it is natural to expect, that one can also derive functional limit theorems under the same assumptions and by using this coupling method.

We next introduce the assumptions on the cone $K$ and on the random walk $\{S(n):n\geq 1\}$. Let $u(x)$ be the
unique strictly positive  on $K$ solution of the following boundary problem:
$$
\Delta u(x)=0,\ x \in K\quad\text{with boundary condition }u\big|_{\partial K}=0.
$$
Let $L_{\mathbb{S}^{d-1}}$ be the Laplace-Beltrami operator on
$\mathbb{S}^{d-1}$ and assume that $\Sigma$ is regular with respect to $L_{\mathbb{S}^{d-1}}$. With this
assumption, there exists 
(see, for example, \cite{BS97})
a complete set of orthonormal eigenfunctions $m_j$  and corresponding
eigenvalues $0<\lambda_1<\lambda_2\le\lambda_3\le\ldots$ satisfying
\begin{align}
 \label{eq.eigen}
  L_{\mathbb{S}^{d-1}}m_j(x)&=-\lambda_jm_j(x),\quad x\in \Sigma\\
  \nonumber m_j(x)&=0, \quad x\in \partial \Sigma .
\end{align}
Define
$$
p:=\sqrt{\lambda_1+(d/2-1)^2}-(d/2-1)>0.
$$
The function $u(x)$ is given by
\begin{equation}
\label{u.from.m}
u(x)=|x|^pm_1\left(\frac{x}{|x|}\right),\quad x\in K.
\end{equation}
As in \cite{DW15}, we shall impose the following conditions on the cone $K$:
\begin{itemize}
\item We assume that there exists an open and connected set
$\widetilde \Sigma\subset \mathbb{S}^{d-1}$ with ${\rm
    dist}(\partial \Sigma, \partial \widetilde \Sigma)>0$ such that
$\Sigma\subset \widetilde \Sigma$ and the
function $m_1$  can be extended to $\widetilde \Sigma$ as a solution to
(\ref{eq.eigen}).
\item $K$ is either convex or starlike (there exists $x_0\in \Sigma$ such that
  $x_0+K\subset K$ and ${\rm dist}(x_0+K, \partial K)>0$). Moreover, K is $C^2$.
\end{itemize}
We impose the following assumptions on the increments of the random walk:
\begin{itemize}
\item {\it Normalisation assumption:} We assume that $\mathbf EX_j=0,
  \mathbf EX_j^2=1,j=1,\ldots,d$. In addition we assume that
  $cov(X_i,X_j)=0$.
\item {\it Moment assumption:} We assume that $\mathbf
  E|X|^{\alpha}<\infty$ with $\alpha=p$ if $p>2$ and some $\alpha>2$
  if $p\le 2$.
\end{itemize}
Let $\{M_K(t),t\in[0,1]\}$ denote the Brownian meander in the cone $K$. Roughly speaking,
this process is the Brownian motion conditioned to stay in $K$ for all $t\in[0,1]$. A rigorous
construction is due to Garbit \cite{Gar09}. Finally, define 
\begin{align*}
K_+:=\{x\in K:&\text{ there exists }\gamma>0\text{ such that for every }R>0\\&\text{ there exists }n \text{ with }\pr(x+S(n)\in D_{R,\gamma},\tau_x>n)>0\}
\end{align*}
and $\tilde{K}_+$ the analogous set where the random walk has step $-X$ instead of $X$. 
\begin{theorem}
\label{t.meander}
For every fixed $x\in K_+$, the process $\{\frac{x+S([nt])}{\sqrt{n}},t\in[0,1]\}$ conditioned
on $\{\tau_x>n\}$ converges towards $\{M_k(t),t\in[0,1]\}$ weakly in
$\left(D[0,1],||\cdot||_{\infty}\right)$.
\end{theorem}
The case of one-dimensional random walks has already been studied in the literature. It is clear
that one has only two non-trivial cones: $(0,\infty)$ and $(-\infty,0)$. Furthermore,
due to the symmetry, it suffices to consider $(0,\infty)$ only. The corresponding invariance
principle for random walks with zero mean and finite variance has been proven by 
Bolthausen \cite{Bolt76}. Doney \cite{Don85} proved a functional limit theorem for asymptotically
stable random walks conditioned to stay positive. 

Shimura \cite{Shim91} has proved convergence in Theorem \ref{t.meander} for two-dimensional walks
with bounded increments. Garbit \cite{Gar11} proved a slight generalisation of Shimura's result.
We are not aware of any further result on convergence towards the Brownian meander in the
multidimensional case.

\vspace{12pt}

We now turn to the weak convergence of random walks conditioned to stay in $K$ at all times.
Such processes are usually defined by using Doob's $h$-transforms. In the case of the standard $d$-dimensional
Brownian motion $B(t)$ one can use the function $u$. Its harmonicity can be rewritten in the following way:
$$
u(x)=\mathbf{E}[u(x+B(t)),\tau_x^{bm}>t],\quad t>0,
$$
where
$$
\tau_x:=\inf\{t>0:x+B(t)\notin K\}.
$$
Then we may consider the probabilistic measure $\mathbf{P}^{(u)}$ given by the following relation:
For any $t>0$ and each continuous and bounded functional $f_t: C[0,t]\mapsto \mathbb{R}$,
$$
\mathbf{E}^{(u)}[f_t(B)|B(0)=x]=\mathbf{E}\left[f_t(B)\frac{u(x+B(t))}{u(x)},\tau^{bm}_x>t\right],\quad x\in K.
$$
As usual, we shall write $\mathbf{P}^{(u)}_x$ for the distribution of the process with starting point $x$.
For our next theorem we need the following property of measures $\mathbf{P}^{(u)}_x$. As $x\to0$,
\begin{equation}
\label{bm.h.conv}
\mathbf{P}^{(u)}_x \text{ converges weakly on }C[0,\infty).
\end{equation}
This convergence is a quite simple consequence of the convergence towards $M_K$
proved by Garbit. In the literature we have found only a one-dimensional version of this convergence: Chaumont \cite{Ch97} has shown this relation for stable processes conditioned
to stay positive. For this reason we give later our proof.

Let $\mathbf{P}^{(u)}_0$ denote the limiting distribution in \eqref{bm.h.conv}. 

In \cite{DW15} a positive harmonic function $V(x)$ for $S(n)$ killed at leaving $K$ has been constructed.
Let $\mathbf{P}^{(V)}$ denote the $h$-transform of $\mathbf{P}$ with the function $V$. 
\begin{theorem}
\label{t.h}
For every fixed $x\in K_+$, the process $\left\{\frac{x+S([nt])}{\sqrt{n}},t\geq0\right\}$ under $\mathbf{P}^{(V)}$
converges weakly to $\{B(t),t\geq0\}$ under $\mathbf{P}^{(u)}_0$. 
\end{theorem}

In the special case of the Weyl chamber of type A this result was proved by Denisov and Wachtel \cite{DW10}. In
the case when $d=1$ and $K=(0,\infty)$ Bryn-Jones and Doney \cite{BD06} proved Theorem \ref{t.h} for random walks
which are integer-valued, aperiodic and belong to the domain of attraction of a standard normal law. Caravenna and Chaumont \cite{CC08} have generalized this result to the
whole class of asymptotically stable random walks.

It is
well-known that for random walks belonging to the domain of attraction of the normal distribution one gets in the limit the three-dimensonal Bessel process. Grabiner
\cite{Grab99} has shown that the radial part of a Brownian motion conditioned to stay in Weyl chambers is also
a Bessel process. It is then immediate from Theorem \ref{t.h} that the rescaled radial part of a random walk converges
to the corresponding Bessel process. This is actually valid in all cones satisfying our geometric assumption.
\begin{corollary}
\label{c.h}
For every fixed $x\in K_+$, the process $\left\{\frac{|x+S([nt])|}{\sqrt{n}},t\geq0\right\}$ under $\mathbf{P}^{(V)}$
converges weakly to a $(2p+d)$-dimensional Bessel process.
\end{corollary}

\vspace{12pt}

We now turn to bridges of random walks conditioned to stay in $K$. 
Here we shall consider lattice random walks. More precisely, we shall assume that $X$ takes values on a
lattice $R$ which is a non-degenerate linear transformation of $\mathbb{Z}^d$. This is weaker than the strong aperiodicity
assumption imposed in \cite{DW15}. This strong aperiodicity was imposed to make direct use of the simplest version of the standard (unconditioned)
local limit theorem from Spitzer's book \cite{Sp76}. But this can be replaced by a local limit theorem proved by Stone
\cite{Stone65}, which is valid for all lattice random walks. Thus, conditioned local limit theorems from \cite{DW15}
remain valid modulo the following changes:
\begin{itemize}
 \item In Theorem 5 one has to take supremum over the set $D_n(x):=\{y\in K: \mathbf{P}(x+S(n)=y)>0\}$.
       In Theorem 6 one has to assume that $y\in D_n(x)$.
 \item If the random walk is not strongly aperiodic then an additional constant depending on the lattice $R$ appears
       in front of the limiting density in (9).      
\end{itemize}

\begin{theorem}
\label{t.bridge}
For all fixed $x\in K_+$ and $y \in \tilde{K}_+$  such that $y\in D_n(x)$, the process $\left\{\frac{x+S([nt])}{\sqrt{n}},t\in[0,1]\right\}$ conditioned on
$\{\tau_x>n,x+S(n)=y\}$ converges weakly to a process $B_K^0(t)$, which we shall call Brownian excursion
in $K$.
\end{theorem}
\begin{remark}
One-dimensional distributions of $B_K^0$ have been calculated in Theorem 6 of \cite{DW15}. The assumption $x\in K_+$ and $y \in \tilde{K}_+$ should be added to the statement of Theorem 6 in \cite{DW15}, because it is necessary for its proof.
\end{remark}

The proof of Theorem \ref{t.bridge} is based on a method suggested by Caravenna and Chaumont \cite{CC13}, who
proved an invariance principle for bridges of one-dimensional random walks conditioned to stay positive. It is worth mentioning that their results are valid for all asymptotically stable random walks. Their
method combines weak convergence towards the meander and local limit theorems for conditioned probabilities. In the case of random walks in cones we have all needed ingredients: Theorem \ref{t.meander} and local limit
theorems from \cite{DW15}.

As usual, invariance principles allow to derive results about weak convergence of appropriate functionals of random processes. If
one knows the limiting behaviour for one particular random walk then the same convergence is valid for all random
walks satisfying the conditions in our theorems. For example, Feierl \cite{Feierl12,Feierl13} has considered
functionals $\max_{k\leq n}(x_d+S_d(k))$ and $\max_{k\leq n}|x_d+S_d(k)-x_1-S_1(k)|$ for random walk bridges in
Weyl chambers $W_A:=\{x\in\mathbb{R}^d:x_1<x_2<\ldots<x_d\}$ and $W_B:=\{x\in\mathbb{R}^d:0<x_1<x_2<\ldots<x_d\}$.
Assuming that every coordinate of $S(n)$ is a simple symmetric random walk, he has shown that the functionals
mentioned above converge weakly after rescaling by $\sqrt{n}$. His proofs are based on very hard calculations,
which can be done for simple random walks only. Our Theorem \ref{t.bridge} implies that one has the same limiting 
distributions for a much larger class of random walks.

\begin{remark}
If the limiting process is a functional of the Brownian motion then one usually proves the weak convergence on the
space of continuous functions. In order to do so, one considers linear interpolations of random walks. For random
walks in cones this is possible only in the case of convex cones, since in a non-convex cone it can happen that a
linear interpolation of points from the cone leaves the cone. This explains the choice of $D$-space with uniform
metric in our theorems. 
\hfill $\diamond$
\end{remark}

\section{Convergence towards the Brownian meander: proof of Theorem \ref{t.meander}}
Let $f:D[0,1]\mapsto\mathbb{R}$ be a non-negative uniformly continuous with respect to the uniform topology function
with values in $[0,1]$. Set also for brevity,
$$
X^{(n)}(t):=\frac{x+S([nt])}{\sqrt{n}},\quad t\geq0.
$$
It suffices to show that
\begin{equation}
\label{t.meander.1}
\mathbf{E}\left[f\left(X^{(n)}\right)\Big|\tau_x>n\right]\to
\mathbf{E}[f(M_K)].
\end{equation}
Similar to \cite{DW15}, we are going to use the strong approximation of random walks by the Brownian motion. We
shall also keep the notation from \cite{DW15}. Let $\varepsilon>0$ be a constant and let
\begin{equation*}
K_{n,\varepsilon}=\{x\in K: {\rm dist}(x,\partial K)\ge n^{1/2-\varepsilon}\}.
\end{equation*}
Define
$$
\nu_n:=\min\{k\geq1: x+S(k)\in K_{n,\varepsilon}\}.
$$
According to Lemma 14 from \cite{DW15},
$$
\mathbf{P}(\nu_n>n^{1-\varepsilon},\tau_x>n^{1-\varepsilon})\leq
\exp\{-Cn^\varepsilon\}.
$$
Furthermore, $\mathbf{P}(\tau_x>n)\sim V(x)n^{-p/2}$ as $n\rightarrow\infty$. Consequently,
\begin{align}
\label{t.meander.2}
\mathbf{E}\left[f\left(X^{(n)}\right),\tau_x>n\right]
=\mathbf{E}\left[f\left(X^{(n)}\right),\nu_n\leq n^{1-\varepsilon},\tau_x>n\right]
+o(\mathbf{P}(\tau_x>n)).
\end{align}

It has been shown by Borovkov \cite{Bor72} that all Fuk-Nagaev inequalities remain
valid for partial maximas of sums of independent random variables. In particular,
the estimates in Corollary 22 from \cite{DW15} hold for $M(n):=\max_{k\leq n}|S(k)|$:
for all $x,y>0$,
\begin{equation}\label{NF3}
\mathbf{P}\left(|M(n)|>x,\max_{k\leq n}|X(k)|\leq y\right)\leq 
2d e^{x/\sqrt{d}y}\left(\frac{\sqrt{d}n}{xy}\right)^{x/\sqrt{d}y}
\end{equation}
and
\begin{equation}\label{NF4}
\mathbf{P}(|M(n)|>x)\leq 
2d e^{x/\sqrt{d}y}\left(\frac{\sqrt{d}n}{xy}\right)^{x/\sqrt{d}y}+n\mathbf{P}(|X(1)|>y).
\end{equation}
Using these bounds in the proof of Lemma 24 from \cite{DW15}, one can easily obtain
\begin{equation}
\label{newbound}
\lim_{n\to\infty}\mathbf{E}\left[|x+S(\nu_n)|^p;\tau_x>\nu_n,
M(\nu_n)>\theta_n\sqrt{n},\nu_n\leq n^{1-\varepsilon}\right] =0.
\end{equation}
Then, by the Markov property,
\begin{align}
\label{t.meander.2a}
\nonumber
&\mathbf{E}\left[f\left(X^{(n)}\right),\nu_n\leq n^{1-\varepsilon},M(\nu_n)>\theta_n\sqrt{n},\tau_x>n\right]\\
\nonumber
&\hspace{1cm}\leq\mathbf{P}\left(\nu_n\leq n^{1-\varepsilon},M(\nu_n)>\theta_n\sqrt{n},\tau_x>n\right)\\
\nonumber
&\hspace{1cm}\leq \frac{C}{n^{p/2}}\mathbf{E}\left[|x+S(\nu_n)|^p;\tau_x>\nu_n,M(\nu_n)>
\theta_n\sqrt{n},\nu_n\leq n^{1-\varepsilon}\right]\\
&\hspace{1cm}=o(n^{-p/2})=o(\mathbf{P}(\tau_x>n)).  
\end{align}

Let
$$
C_n:=\{\nu_n\leq n^{1-\epsilon}, M(\nu_n)\leq \theta_n\sqrt{n}\}.
$$
Define the functions
$$
f(y,k,X^{(n)}) = f\left(\frac{y}{\sqrt{n}}\textbf{1}_{\{t\leq \frac{k}{n}\}}+
X^{(n)}(t)\textbf{1}_{\{t> \frac{k}{n}\}}\right)
$$
Then it holds 
\begin{align*}
&\sup_{t\in[0,1]}\left|\left(\frac{x+S(\nu_n)}{\sqrt{n}}\textbf{1}_{\{t\leq \frac{k}{n}\}}
+X^{(n)}(t))\textbf{1}_{\{t>\frac{k}{n}\}})\right)-X_n(t)\right|\\
&\hspace{2cm} \leq \frac{\max_{k\leq \nu_n}|S(\nu_n)-S(k)|}{\sqrt{n}}\leq \frac{2M(\nu_n)}{\sqrt{n}}.
\end{align*}
On the set $\{M(\nu_n)\leq \theta_n\sqrt{n}\}$ this is smaller than $2\theta_n$. Then, as $n\to\infty$, 
$$
|f(X^{(n)}-f(S(\nu_n),\nu_n,X^{(n)})| = o(1),\text{ uniformly on }C_n,
$$
since $\theta_n\to 0, n\to\infty$.
From this estimate we infer that
\begin{align*}
&\mathbf{E}\left[f\left(X^{(n)}\right),C_n,\tau_x>n\right]\\
&\hspace{1cm}=\mathbf{E}\left[f\left(x+S(\nu_n),\nu_n,X^{(n)}\right),C_n,\tau_x>n\right]
+o(\mathbf{P}(\tau_x>n)).
\end{align*}
Then, taking into account \eqref{t.meander.2} and \eqref{t.meander.2a}, we have
\begin{align}
\label{t.meander.3}
\nonumber
&\mathbf{E}\left[f\left(X^{(n)}\right),\tau_x>n\right]\\
&\hspace{1cm}=\mathbf{E}\left[f\left(x+S(\nu_n),\nu_n,X^{(n)}\right),C_n,\tau_x>n\right]
+o(\mathbf{P}(\tau_x>n)).
\end{align}
Thus, to prove the theorem, it suffices to consider the right hand side of the above equation. 
Note first that we can write 
\begin{align*}
&\mathbf{E}\left[f\left(\frac{x+S(\nu_n)}{\sqrt{n}},\nu_n, X^{(n)}\right), C_n,\tau_x>n\right]\\
&=\sum_{k\leq n^{1-\varepsilon}}\int_{K_{n,\varepsilon}}
\mathbf{P}(\nu_n = k,\tau_x>k, M(k)\leq \theta_n\sqrt{n},x+S(k)\in dz)\mathbf{E}[f(z,k,X^{(n)}),\tau_z>n-k]. 
\end{align*}
We now note that it is sufficient to show that, uniformly in $z\in K_{n,\varepsilon}$, $k\leq n^{1-\varepsilon}$,
$|z|\leq \theta_n\sqrt{n}$,
\begin{equation}
\label{t.meander.4}
\mathbf{E}\left[f\left(z,k,X^{(n)}\right),\tau_z>n-k\right]
=(1+o(1))\mathbf{E}[f(M_K)]\frac{u(z)}{n^{p/2}}.
\end{equation}
Indeed, \eqref{t.meander.4} implies that
\begin{align*}
&\mathbf{E}[f(x+S(\nu_n),\nu_n, X^{(n)}), C_n,\tau_x>n]\\
&\hspace{1cm}\sim\mathbf{E}[f(M_K)]
\frac{1}{n^{p/2}}\mathbf{E}\left[u(x+S(\nu_n)), C_n,\tau_x>\nu_n\right].
\end{align*}
It follows from \eqref{newbound} and Lemma 21 in \cite{DW15} that
\begin{equation}
\label{t.meander.4a}
\mathbf{E}\left[u(x+S(\nu_n)), C_n,\tau_x>\nu_n\right]
\sim \mathbf{E}\left[u(x+S(\nu_n)),\tau_x>\nu_n\right]\sim V(x).
\end{equation}
Now we infer from \eqref{t.meander.3} that
$$
\mathbf{E}\left[f\left(X^{(n)}\right),\tau_x>n\right]
\sim \mathbf{E}[f(M_K)]\frac{V(x)}{n^{p/2}}
\sim \mathbf{E}[f(M_K)]\mathbf{P}(\tau_x>n).
$$
In other words,
$$
\mathbf{E}\left[f\left(X^{(n)}\right)\Big|\tau_x>n\right]\sim \mathbf{E}[f(M_K)].
$$
In order to prove \eqref{t.meander.4} we first note that, according to Lemma 17 from \cite{DW15}, one can define
on a joint probability space a copy of $S(n)$ and a Brownian motion $B(t)$ on the same space, such that if 
$\mathbf{E}|X|^{2+\delta}$ is finite then for any $\gamma$ so that  $0<\gamma<\frac{\delta}{2(2+\delta)}$ the probability of the event
$$
\left\{\sup_{u\leq n}|S([u])-B(u))|> n^{1/2-\gamma}\right\}
$$
does not exceed $n^{-r}$ with $r=r(\delta,\gamma)=2\gamma+\gamma\delta-\delta/2$. Setting
$$
A_n:=\left\{\sup_{u\leq n}|S([u])-B(u))|\leq n^{1/2-\gamma}\right\},
$$
we have
$$
\mathbf{E}\left[f\left(z,k,X^{(n)}\right),\tau_z>n-k\right]=
\mathbf{E}\left[f\left(z,k,X^{(n)}\right),A_n,\tau_z>n-k\right]+O(n^{-r}).
$$
In view of the proof of Lemma 20 in \cite{DW15}, for all $\varepsilon>0$ sufficiently small,
$$
\frac{n^{-r}}{\mathbf{P}(\tau_z>n)}\to0
$$
uniformly in $z\in K_{n,\varepsilon}$. Consequently,
\begin{align}
\label{t.meander.5}
\nonumber
&\mathbf{E}\left[f\left(z,k,X^{(n)}\right),\tau_z>n-k\right]\\
&\hspace{1cm}=\mathbf{E}\left[f\left(z,k,X^{(n)}\right),A_n,\tau_z>n-k\right]+o(\mathbf{P}(\tau_z>n)).
\end{align}
For every $z\in K_{n,\varepsilon}$ we define
$$
z^\pm=z\pm R_0x_0 n^{1/2-\gamma},
$$
where $x_0$ is such that $|x_0|=1$, $x_0+K\subset K$ and $R_0$ is such that ${\rm dist}(R_0x_0+K,\partial K)>1$.
Note also that this choice of $R_0$ ensures that $R_0x_0n^{1/2-\gamma}\subset K_{n,\gamma}$.

If we take $\varepsilon<\gamma$, then for any $\varepsilon'>\varepsilon$ there exists a positive integer $n(\varepsilon')$ such that
$z^{\pm}\in K_{n,\varepsilon'}$ as soon as $n\geq n(\varepsilon')$ and $z\in K_{n,\varepsilon}$.

Let $B^{(n)}(t):=n^{-1/2}B(nt)$, $t\geq0$. From the coupling described above we have the following relations
$$
\{\tau^{bm}_{z^-}>n-k\}\cap A_n\subset \{\tau_{z}>n-k\}\cap A_n\subset \{\tau^{bm}_{z^+}>n-k\}\cap A_n
$$
and
$$
\left|\widetilde{f}\left(z^\pm,k,B^{(n)}\right)-f\left(z,k,X^{(n)}\right)\right|=o(1)
\quad\text{uniformly on }A_n,
$$
where
$$
\widetilde{f}\left(z,k,B^{(n)}\right):=
f\left(\frac{z}{\sqrt{n}}+B^{(n)}(t-k/n){\bf 1}\{t>k/n\}\right)
$$
Moreover, the latter relation is uniform in $k$, $z\in K_{n,\varepsilon}$ and in $A_n$.
Combining these relations with \eqref{t.meander.5}, we obtain
\begin{align}
\label{t.meander.6}
\nonumber
&\mathbf{E}\left[\widetilde{f}\left(z^-,k,B^{(n)}\right),\tau^{bm}_{z^-}>n\right]+o(\mathbf{P}(\tau_z>n))\\
&\hspace{1cm}\leq \mathbf{E}\left[f\left(z,k,X^{(n)}\right),\tau_z>n-k\right]+o(\mathbf{P}(\tau_z>n))\\
\nonumber
&\hspace{2cm}\leq\mathbf{E}\left[\widetilde{f}\left(z^+,k,B^{(n)}\right),\tau^{bm}_{z^+}>n-n^{1-\varepsilon}\right]
+o(\mathbf{P}(\tau_z>n)),
\end{align}
uniformly in $k$, $z\in K_{n,\varepsilon}$.
Moreover, it follows from the main result in Garbit \cite{Gar09} that for every continuous functional $g$
$$
\sup_{y\in K:|y|\leq r}\mathbf{E}[g(y+B)|\tau^{bm}_y>1]\to\mathbf{E}[g(M_K)],\quad\text{as }r\to0.
$$
Then, by the scaling property of the Brownian motion and the definition of $\widetilde{f}(y,k,\cdot)$,
$$
\max_{k\leq n^{1-\varepsilon}}\sup_{z\in K_{n,\varepsilon}:|z|\leq\theta_n\sqrt{n}}
\Big|\mathbf{E}\left[\widetilde{f}\left(z^+,k,B^{(n)}\right)\big|\tau^{bm}_{z^+}>n-n^{1-\varepsilon}\right]
-\mathbf{E}[f(M_K)]\Big|\to0
$$
and
$$
\max_{k\leq n^{1-\varepsilon}}\sup_{z\in K_{n,\varepsilon}:|z|\leq\theta_n\sqrt{n}}
\Big|\mathbf{E}\left[\widetilde{f}\left(z^-,k,B^{(n)}\right)\big|\tau^{bm}_{z^-}>n\right]-\mathbf{E}[f(M_K)]\Big|\to0.
$$
According to Lemmata 18 and 20 from \cite{DW15},
$$
\mathbf{P}(\tau^{bm}_{z^+}>n-n^{1-\varepsilon})\sim\mathbf{P}(\tau^{bm}_{z^-}>n)\sim\varkappa\frac{u(z)}{n^{p/2}}
$$
and
$$
\mathbf{P}(\tau_{z}>n)\sim\varkappa\frac{u(z)}{n^{p/2}}
$$
uniformly in $z\in K_{n,\varepsilon}$.

Combining these relations with \eqref{t.meander.6}, we arrive at \eqref{t.meander.4}. This completes the proof of the theorem.

\section{Convergence of $h$-transforms: proof of Theorem \ref{t.h}}
\subsection{Proof of \eqref{bm.h.conv}}
Due to the scaling property of the Brownian motion, it suffices to prove the convergence on $C[0,1]$. Let $f$ be a
bounded and uniformly continuous function on $C[0,1]$ with values in $[0,1]$. Fix some $R>1$ and split
\begin{equation}
\label{proof_of_3.1}
\mathbf{E}_x^{(u)}[f(B)]=\mathbf{E}_x^{(u)}[f(B),|x+B(1)|\leq R]
+\mathbf{E}_x^{(u)}[f(B),|x+B(1)|>R].
\end{equation}
In order to bound the second summand we note that we have good control over the density $\frac{\mathbf{P}(x+B(1)\in dz,\tau^{bm}_x>1)}{dz}=:p(1,x,z)$. Namely, it
follows from the proof of Lemma 5.4 in \cite{Gar09} that
\begin{equation}\label{eq:garbit-uniform}
\sup_{x\in K, |x|\leq\frac{1}{2}}\frac{p(1,x,z)}{u(x)}\leq Ce^{-\gamma|z|^2}
\end{equation}
for some suitable $C,\gamma>0$. Consequently, introducing spherical coordinates and using the estimate $u(x)\leq C|x|^p$, we
obtain
\begin{align*}
\mathbf{E}_x^{(u)}[f(B),|x+B(1)|>R]
&=\frac{1}{u(x)}\mathbf{E}[f(x+B)u(x+B(1)),\tau^{bm}_x>1,|x+B(1)|>R]\\
&\leq\frac{C}{u(x)}\mathbf{E}[|x+B(1)|^p,\tau^{bm}_x>1,|x+B(1)|>R]\\
&\leq C\int_R^\infty r^{p+d-1}e^{-\gamma r^2}dr.
\end{align*}
uniformly in $x\in K, |x|\leq \frac{1}{2}$. 
Therefore,
\begin{equation}
\label{proof_of_3.2}
\sup_{x\in K, |x|\leq\frac{1}{2}}\mathbf{E}_x^{(u)}[f(B),|x+B(1)|>R]\to0
\quad\text{as }R\to\infty.
\end{equation}
For the first summand in \eqref{proof_of_3.1} we have the following representation
\begin{align*}
&\mathbf{E}_x^{(u)}[f(B),|x+B(1)|\leq R]\\
&\hspace{1cm}=\frac{\mathbf{P}(\tau^{bm}_x>1)}{u(x)}
\mathbf{E}[f(x+B)u(x+B(1)){\rm 1}\{|x+B(1)|\leq R\}|\tau^{bm}_x>1].
\end{align*}
The functional $f(x+B)u(x+B(1)){\rm 1}\{|x+B(1)|\leq R\}$ is bounded and the set of its discontinuities is a null set with respect to the distribution of the meander $M_K$.
Thus, the convergence result from \cite{Gar09} is still applicable and, consequently,
\begin{align*}
&\mathbf{E}[f(x+B)u(x+B(1)){\rm 1}\{|x+B(1)|\leq R\}|\tau^{bm}_x>1]\\
&\hspace{1cm}\to
\mathbf{E}[f(M_K)u(M_K(1)){\rm 1}\{M_K(1)\leq R\}],\quad x\to0.
\end{align*}
Recalling that $\mathbf{P}(\tau^{bm}_x>1)\sim \varkappa u(x)$, we then have
$$
\mathbf{E}_x^{(u)}[f(B),|x+B(1)|\leq R]\to 
\varkappa \mathbf{E}[f(M_K)u(M_K(1)){\rm 1}\{M_K(1)\leq R\}]
$$
Combining this with \eqref{proof_of_3.2} and noting that, by monotone convergence,
$$
\mathbf{E}[f(M_K)u(M_K(1)){\rm 1}\{M_K(1)\leq R\}]\to
\mathbf{E}[f(M_K)u(M_K(1))]\quad\text{as }R\to\infty,
$$
we finally get
\begin{equation}
\label{h-trans.final}
\lim_{|x|\to0}\mathbf{E}_x^{(u)}[f(B)]=\varkappa\mathbf{E}[f(M_K)u(M_K(1))].
\end{equation}
\subsection{Proof of Theorem \ref{t.h}}
Take a uniformly continuous and bounded functional in $\left(D[0,1],||\cdot||_{\infty}\right)$ with $0\leq f\leq 1$.
We want to show that
$$
\mathbf{E}_x^{(V)}[f(X^{(n)})]\to\mathbf{E}^{(u)}_0[f(B)] \quad\text{as } n\to\infty.
$$
We first note that 
$$
\mathbf{E}_x^{(V)}[f(X^{(n)}),\nu_n>n^{1-\varepsilon}]
\leq \frac{1}{V(x)}\mathbf{E}[V(x+S(n)),\tau_x>n,\nu_n>n^{1-\varepsilon}].
$$
By the Markov property and the harmonicity of $V$ we get
\begin{align*}
&\mathbf{E}[V(x+S(n)),\tau_x>n,\nu_n>n^{1-\epsilon}]\\
&\hspace{1cm}=
\mathbf{E}[\mathbf{E}[V(x+S(n))\textbf{1}\{\tau_x>n\}|\mathcal{F}_{n^{1-\varepsilon}}],\nu_n>n^{1-\varepsilon}]\\
&\hspace{1cm}= \mathbf{E}[V(x+S(n^{1-\varepsilon})),\tau_x>n^{1-\varepsilon},\nu_n>n^{1-\varepsilon}].
\end{align*}
We recall the definition of $V$ as
\[
V(x) = u(x) -\mathbf{E}[u(x+S(\tau_x))]+\mathbf{E}[\sum_{k=0}^{\tau_x-1}f(x+S(k))].
\]
From estimates found in the proof of Lemma 8 and 12 in \cite{DW15} one gets the estimate
\begin{equation}\label{eq:estimateV}
V(x)\leq u(x)+C|x|^{p-\delta},\quad x\in K, |x|\geq 1
\end{equation}
for some $\delta>0$ small enough. By Lemma 16 of \cite{DW15} we have 
$$
\mathbf{E}[u(x+S(n^{1-\varepsilon})),\tau_x>n^{1-\varepsilon},\nu_n>n^{1-\epsilon}]
\leq C(x)e^{-cn^{\varepsilon}}.
$$
Furthermore, applying H\"older inequality we have 
\begin{align*}
&\mathbf{E}[|x+S(n^{1-\varepsilon})|^{p-\delta},\tau_x>n^{1-\varepsilon},\nu_n>n^{1-\epsilon}]\\
&\hspace{1cm}\leq \left(\mathbf{E}[|x+S(n^{1-\varepsilon})|^{p}]\right)^{\frac{p-\delta}{p}}
\left(\mathbf{P}(\tau_x>n^{1-\varepsilon},\nu_n>n^{1-\varepsilon})\right)^{\frac{\delta}{p}}.
\end{align*}
Taking now into account Lemma 14 in \cite{DW15} we get for suitably small $\varepsilon>0$ 
\begin{align*}
&\mathbf{E}[|x+S(n^{1-\varepsilon})|^{p-\delta},\tau_x>n^{1-\varepsilon},\nu_n>n^{1-\epsilon}]\\
&\hspace{1cm}\leq C(x)e^{-cn^{\varepsilon}}.
\end{align*}
As a result we get
\begin{equation}
\label{eq:partial1}
\mathbf{E}_x^{(V)}[f(X^{(n)}),\nu_n>n^{1-\varepsilon}]=o(1),\quad n\to\infty.
\end{equation}
Using the Markov property and the harmonicity of $V$ once again, we get 
\begin{align*}
&\mathbf{E}_x^{(V)}[f(X^{(n)}),\nu_n\leq n^{1-\varepsilon}, M(\nu_n)>\theta_n\sqrt{n}]\\
&\hspace{1cm}\leq
\frac{1}{V(x)}\mathbf{E}[V(x+S(n)),\nu_n\leq n^{1-\varepsilon},M(\nu_n)>\theta_n\sqrt{n},\tau_x>n]\\
&\hspace{1cm}\leq
\frac{1}{V(x)}\mathbf{E}[V(x+S(\nu_n)),\nu_n\leq n^{1-\varepsilon},M(\nu_n)>\theta_n\sqrt{n},\tau_x>\nu_n].
\end{align*}
It is immediate from \eqref{eq:estimateV} that 
$$
V(x)\leq c|x|^p,\quad x\in K, |x|\geq 1
$$
Consequently, 
\begin{align*}
&\mathbf{E}_x^{(V)}[f(X^{(n)}),\nu_n\leq n^{1-\varepsilon}, M(\nu_n)>\theta_n\sqrt{n}]\\
&\hspace{1cm}\leq 
\frac{C}{V(x)}\mathbf{E}[|x+S(\nu_n)|^p,\nu_n\leq n^{1-\varepsilon},M(\nu_n)>\theta_n\sqrt{n},\tau_x>\nu_n].
\end{align*}
Then, in view of \eqref{newbound},
\begin{equation}\label{eq:partial2}
\mathbf{E}_x^{(V)}[f(X^{(n)}),\nu_n\leq n^{1-\varepsilon}, M(\nu_n)>\theta_n\sqrt{n}] = o(1),\quad n\to \infty.
\end{equation}

Together, \eqref{eq:partial1} and \eqref{eq:partial2} imply
$$
\mathbf{E}_x^{(V)}[f(X^{(n)})]=
\mathbf{E}_x^{(V)}[f(X^{(n)}),\nu_n\leq n^{1-\varepsilon}, M(\nu_n)\leq \theta_n\sqrt{n}] + o(1),
\quad n\to \infty.
$$
Therefore, we have to look only at convergence of 
$$
\mathbf{E}_x^{(V)}[f(X^{(n)}),\nu_n\leq n^{1-\varepsilon}, M(\nu_n)\leq \theta_n\sqrt{n}].
$$
We have already seen in the proof of Theorem \ref{t.meander} that 
$$
|f(X^{(n)})-f(x+S(\nu_n),\nu_n,X^{(n)})| = o(1)\text{ uniformly on }
\{\nu_n\leq n^{1-\varepsilon},M(\nu_n)\leq \theta_n\sqrt{n}\}.
$$
Therefore, 
\begin{align*}
&\mathbf{E}_x^{(V)}[f(X^{(n)}),\nu_n\leq n^{1-\varepsilon}, M(\nu_n)\leq \theta_n\sqrt{n}]\\
&\hspace{1cm}= 
\mathbf{E}_x^{(V)}[f(x+S(\nu_n),\nu_n,X^{(n)}),\nu_n\leq n^{1-\varepsilon}, M(\nu_n)\leq \theta_n\sqrt{n}] + o(1).
\end{align*}
We can then continue with 
\begin{align}
\label{t.h.1a}
&\mathbf{E}_x^{(V)}[f(x+S(\nu_n),\nu_n,X^{(n)}),\nu_n\leq n^{1-\varepsilon}, M(\nu_n)\leq \theta_n\sqrt{n}]\\
\nonumber
&\hspace{1cm}=\sum_{k=1}^{n^{1-\varepsilon}}\frac{1}{V(x)}\int_{K_{n,\varepsilon}} 
\mathbf{P}(x+S(k)\in dz, \tau_x>k,\nu_n=k, M(k)\leq\theta_n\sqrt{n})\\
\nonumber
&\hspace{4.2cm}\times\mathbf{E}[f(z,k,X^{(n)})V(z+S(n-k)),\tau_z>n-k].
\end{align}
To estimate $\mathbf{E}[f(z,k,X^{(n)})V(z+S(n-k)),\tau_z>n-k]$ for $k\leq n^{1-\epsilon}$ we shall use once again
coupling arguments. We first show that 
\begin{equation}
\label{t.h.2}
\mathbf{E}[V(z+S(n-k)),\tau_z>n-k,A_n^c] = o(u(z))
\end{equation}
and
\begin{equation}
\label{t.h.2a}
\mathbf{E}[u(z+B(n-k)),\tau^{bm}_z>n-k,A_n^c] = o(u(z))
\end{equation}
uniformly for $k\leq n^{1-\varepsilon}, z\in K_{n,\varepsilon}, |z|\leq \theta_n\sqrt{n}.$ Fix some $\eta>0$ and
split
\begin{align*}
&\mathbf{E}[V(z+S(n-k)),\tau_y>n-k,A_n^c]\\
&\hspace{1cm}=\mathbf{E}[V(z+S(n-k)),\tau_z>n-k,|z+S(n-k)|\leq n^{\frac{1}{2}+\eta},A_n^c]\\
&\hspace{3cm}+\mathbf{E}[V(z+S(n-k)),\tau_z>n-k,|z+S(n-k)|> n^{\frac{1}{2}+\eta},A_n^c]\\
&\hspace{1cm}=: E_{n,1}+E_{n,2}.
\end{align*}
Note that
$$
E_{n,1}\leq Cn^{\frac{p}{2}+p\eta}\mathbf{P}(A_n^c)\leq Cn^{\frac{p}{2}+p\eta-r}.
$$
It follows easily from Lemma 19 in \cite{DW15} that
$$
u(z)\geq C n^{p/2-p\varepsilon},\ z\in K_{n,\varepsilon}.
$$
Then, choosing $\varepsilon$ and $\eta$ sufficiently small, we get 
$$
E_{n,1}=o(u(z)).
$$
Further,
\begin{align*}
E_{n,2}&\leq \mathbf{E}[V(z+S(n-k)),|z+S(n-k)|> n^{\frac{1}{2}+\eta}]\\
&\leq C\mathbf{E}[(M(n))^p, M(n)> n^{\frac{1}{2}+\eta}].
\end{align*}
Integrating by parts, we get
\begin{align*}
&\mathbf{E}[(M(n))^p, M(n)> n^{\frac{1}{2}+\eta}]\\
&\hspace{1cm}=n^{p(1/2+\eta)}\mathbf{P}(M(n)>n^{1/2+\eta})
+p\int_{n^{1/2+\eta}}^\infty x^{p-1}\mathbf{P}(M(n)>x)dx.
\end{align*}
Using now \eqref{NF4} with $y=\frac{\varepsilon}{\sqrt{d}}x$ and taking into account the moment assumption,
one can easily get
$$
\mathbf{E}[(M(n))^p, M(n)> n^{\frac{1}{2}+\eta}]=o(n^{p/2-p\varepsilon})=o(u(z)).
$$
This finishes the proof of \eqref{t.h.2}. The proof of \eqref{t.h.2a} is even simpler.

As we have seen in the proof of Theorem \ref{t.meander},
$$
\{\tau_{z^-}^{bm}>n\}\cap A_n \subset\{\tau_{z}>n\}\cap A_n\subset \{\tau_{z^+}^{bm}>n\}\cap A_n.
$$
Furthermore, we have the estimate 
\begin{equation}\label{eq:longestimateV}
|V(x)-u(x)|\leq C(1+|x|^{p-\delta}),\quad x\in K
\end{equation}
which follows from steps from the proofs of Lemma 8 and 12 in \cite{DW15}. Therefore,
\begin{align*}
&\mathbf{E}[f(z,k,X^{(n)})V(z+S(n-k)),\tau_z>n-k,A_n]\\
&\hspace{1cm}\geq \mathbf{E}[f(z,k,X^{(n)})u(z+S(n-k)),\tau_{z^-}^{bm}>n-k,A_n]\\
&\hspace{2cm}-C\mathbf{E}[(1+|z+S(n-k)|^{p-\delta}),\tau^{bm}_{z^-}>n-k,A_n].
\end{align*}
Note now that
$$
|(z+S(n-k))-(z^-+B(n-k))|^{p-\delta}\leq n^{(1/2-\gamma)(p-\delta)} 
$$
on the event $A_n$. As a result,
\begin{align*}
&\mathbf{E}[V(z+S(n-k))f(z,k,X^{(n)}),\tau_z>n-k,A_n]\\
&\hspace{1cm}\geq \mathbf{E}[f(z,k,X^{(n)})u(z+S(n-k)),\tau_{z^-}^{bm}>n-k,A_n]\\
&\hspace{2cm}-C\mathbf{E}[|z^-+B(n-k)|^{p-\delta},\tau^{bm}_{z^-}>n-k]\\
&\hspace{2cm}-Cn^{(1/2-\gamma)(p-\delta)}\mathbf{P}(\tau^{bm}_{z^-}>n-k).
\end{align*}
Using the scaling property of the Brownian motion and applying \eqref{eq:garbit-uniform},
we get, uniformly in $k\leq n^{1-\varepsilon}$ and $|y|=o(\sqrt{n})$ for $q>0$,
\begin{align}
\label{t.h.3}
\nonumber
\mathbf{E}[(y+B(n-k))^q,\tau^{bm}_y>n-k]
&=(n-k)^{q/2}\mathbf{E}\left[\left(\frac{y}{\sqrt{n-k}}+B(1)\right)^q,\tau^{bm}_{\frac{y}{\sqrt{n-k}}}>1\right]\\
\nonumber 
&\leq C(n-k)^{q/2}u\left(\frac{y}{\sqrt{n-k}}\right)\\
&\leq Cn^{(q-p)/2}u(y),\quad q>0.
\end{align}
Using this inequality with $q=p-\delta$ and recalling that $\mathbf{P}(\tau^{bm}_{z^-}>n-k)\sim u(z^-)n^{-p/2}$,
we have
\begin{align*}
&\mathbf{E}[V(z+S(n-k))f(z,k,X^{(n)}),\tau_z>n-k,A_n]\\
&\hspace{1cm}\geq \mathbf{E}[f(z,k,X^{(n)})u(z+S(n-k)),\tau_{z^-}^{bm}>n-k,A_n]+o(u(z^-))
\end{align*}
uniformly in $k\leq n^{1-\varepsilon}$ and $z\in K_{n,\varepsilon}$ with $|z|\leq\theta_n\sqrt{n}$.
Noting that (43) in \cite{DW15} implies that $u(z^\pm)\sim u(z)$ for all $z\in K_{n,\varepsilon}$ with
$|z|\leq\sqrt{n}$, we finally get
\begin{align*}
&\mathbf{E}[V(z+S(n-k))f(z,k,X^{(n)}),\tau_z>n-k,A_n]\\
&\hspace{1cm}\geq \mathbf{E}[f(z,k,X^{(n)})u(z+S(n-k)),\tau_{z^-}^{bm}>n-k,A_n]+o(u(z^-)).
\end{align*}

We next apply the relation
\begin{align}
\label{t.h.4}
\nonumber
&\Big|\mathbf{E}[u(z+S(n-k)),\tau_{z^-}^{bm}>n-k,A_n]\\
&\hspace{1cm} -\mathbf{E}[u(z^-+B(n-k)),\tau_{z^-}^{bm}>n-k,A_n]\Big|=o(u(z)),
\end{align}
which will be proved later. This estimate yields
\begin{align*}
&\mathbf{E}[V(z+S(n-k))f(z,k,X^{(n)}),\tau_z>n-k,A_n]\\
&\hspace{1cm}\geq \mathbf{E}[f(z,k,X^{(n)})u(z^-+B(n-k)),\tau_{z^-}^{bm}>n-k,A_n]+o(u(z))
\end{align*}
uniformly in $k\leq n^{1-\varepsilon}$ and $z\in K_{n,\varepsilon}$ with $|z|\leq\theta_n\sqrt{n}$.

As we have seen in the proof of Theorem \ref{t.meander},
$$
f(z,k,X^{(n)})-\widetilde{f}\left(z^-,k,B^{(n)}\right)=o(1)
$$
uniformly in $k\leq n^{1-\varepsilon}$, $z\in K_{n,\varepsilon}$ with $|z|\leq\theta_n\sqrt{n}$
and in $A_n$. This implies that
\begin{align*}
&\mathbf{E}[f(z,k,X^{(n)})u(z^-+B(n-k)),\tau_{z^-}^{bm}>n-k,A_n]\\
&\hspace{2cm}-\mathbf{E}[\widetilde{f}\left(z^-,k,B^{(n)}\right)u(z^-+B(n-k)),\tau_{z^-}^{bm}>n-k,A_n]\\
&\hspace{1cm}=o\left(\mathbf{E}[u(z^-+B(n-k)),\tau_{z^-}^{bm}>n-k]\right)=o(u(z)).
\end{align*}
In the last step we used the harmonicity of $u$ and the relation $u(z^-)\sim u(z)$. From this estimate
and \eqref{t.h.2a} we infer
\begin{align*}
&\mathbf{E}[V(z+S(n-k))f(z,k,X^{(n)}),\tau_z>n-k,A_n]\\
&\hspace{1cm}\geq \mathbf{E}\left[\widetilde{f}\left(z^-,k,B^{(n)}\right)u(z^-+B(n-k)),\tau_{z^-}^{bm}>n-k,A_n\right]+o(u(z))\\
&\hspace{1cm}=\mathbf{E}\left[\widetilde{f}\left(z^-,k,B^{(n)}\right)u(z^-+B(n-k)),\tau_{z^-}^{bm}>n-k\right]+o(u(z))\\
&\hspace{1cm}=u(z)\mathbf{E}^{(u)}_{z^-/\sqrt{n}}\left[\widetilde{f}\left(z^-,k,B^{(n)}\right)\right]+o(u(z))
\end{align*}
Then, in view of convergence \eqref{bm.h.conv},
\begin{align*}
\mathbf{E}[V(z+S(n-k))f(z,k,X^{(n)}),\tau_z>n-k,A_n]\geq u(z)\mathbf{E}_0^{(u)}[f(B)]+o(u(z)),
\end{align*}
uniformly in $k\leq n^{1-\varepsilon}$, $z\in K_{n,\varepsilon}$ with $|z|\leq\theta_n\sqrt{n}$.
By similar arguments,
\begin{align*}
\mathbf{E}[V(z+S(n-k))f(z,k,X^{(n)}),\tau_z>n-k,A_n]\leq u(z)\mathbf{E}_0^{(u)}[f(B)]+o(u(z)),
\end{align*}
uniformly in $k\leq n^{1-\varepsilon}$, $z\in K_{n,\varepsilon}$ with $|z|\leq\theta_n\sqrt{n}$.

Combining these inequalities and \eqref{t.h.1a}, we obtain
\begin{align*}
&\mathbf{E}_x^{(V)}[f(S(\nu_n),\nu_n,X^{(n)}),\nu_n\leq n^{1-\varepsilon}, M(\nu_n)\leq \theta_n\sqrt{n}]\\
&\hspace{1cm}=\frac{\mathbf{E}^{(u)}_0[f(B)]+o(1)}{V(x)}\mathbf{E}\left[u(x+S(\nu_n)),
\nu_n\leq n^{1-\varepsilon}, M(\nu_n)\leq \theta_n\sqrt{n},\tau_x>\nu_n\right].
\end{align*}
Taking into account \eqref{t.meander.4a}, we have
\begin{align*}
\mathbf{E}_x^{(V)}[f(S(\nu_n),\nu_n,X^{(n)}),\nu_n\leq n^{1-\varepsilon}, M(\nu_n)\leq \theta_n\sqrt{n}]
=\mathbf{E}^{(u)}_0[f(B)]+o(1).
\end{align*}
Combining this with \eqref{eq:partial1} and \eqref{eq:partial2}, we finally obtain
$$
\mathbf{E}_x^{(V)}[f(X^{(n)})]\to \mathbf{E}^{(u)}_0[f(B)].
$$
It remains to show \eqref{t.h.4}. 
It follows from Lemma 7 in \cite{DW15} that
\begin{equation}
\label{t.h.5}
|\nabla u(x)|\leq C|x|^{p-1},\ x\in K.
\end{equation}

Then, by the Taylor formula, 
\begin{align*}
&|u(z+S(n-k))-u(z^-+B(n-k))|\\
&\leq C\bigl(|z^--z|+|B(n-k)-S(n-k)|\bigr)\left(|z+S(n-k)|^{p-1}+|z^-+B(n-k)|^{p-1}\right).
\end{align*}
From definitions of $z^-$ and $A_n$, we get for $p>1$ the bound
\begin{align*}
&|u(z+S(n-k))-u(z^-+B(n-k))|\\
&\hspace{1cm}\leq C n^{1/2-\gamma}\left(|z+S(n-k)|^{p-1}+|z^-+B(n-k)|^{p-1}\right)\\
&\hspace{1cm}\leq Cn^{p(1/2-\gamma)}+Cn^{1/2-\gamma}|z^-+B(n-k)|^{p-1}.
\end{align*}
Then
\begin{align*}
&\Big|\mathbf{E}[u(z+S(n-k)),\tau_{z^-}^{bm}>n-k,A_n]\\
&\hspace{2cm} -\mathbf{E}[u(z^-+B(n-k)),\tau_{z^-}^{bm}>n-k,A_n]\Big|\\
&\hspace{1cm}\leq Cn^{p(1/2-\gamma)}+Cn^{1/2-\gamma}\mathbf{E}[|z^-+B(n-k)|^{p-1},\tau^{bm}_{z^-}>n-k].
\end{align*}
Using now \eqref{t.h.3} with $q=p-1$, we coclude that \eqref{t.h.4} is proved for $p\geq1$. 

Consider now the case $p<1$.
If $|z^-+B(n-k)|\leq \theta_n\sqrt{n}$ with some $\theta_n\to0$ then
\begin{align}
\label{t.h.6}
\nonumber
&\mathbf{E}[u(z^-+B(n-k)),\tau^{bm}_{z^-}>n-k,|z^-+B(n-k)|\leq \theta_n\sqrt{n}]\\
&\hspace{1cm}\leq \theta_n^p n^{p/2}\mathbf{P}(\tau^{bm}_{z^-}>n-k)=o(u(z)).
\end{align}
Moreover, $|z^-+B(n-k)|\leq \theta_n\sqrt{n}$ implies that
$$
|z+S(n-k)|\leq 2\theta_n\sqrt{n}\quad \text{on }A_n
$$
for any sequence $\theta_n$, which goes to zero sufficiently slow. Consequently,
\begin{equation}
\label{t.h.7}
\mathbf{E}[u(z+S(n-k)),\tau_{z^-}^{bm}>n-k,A_n,|z^-+B(n-k)|\leq \theta_n\sqrt{n}]=o(u(z)).
\end{equation}
If $|z^-+B(n-k)|> \theta_n\sqrt{n}$ then
$$
|z+S(n-k)|\geq \frac{1}{2}\theta_n\sqrt{n}\quad \text{on }A_n,
$$
for any sequence $\theta_n$, which goes to zero sufficiently slow.
Then it follows from the Taylor formula and \eqref{t.h.5} that if $p<1$ then
$$
|u(z+S(n-k))-u(z^-+B(n-k))|\leq Cn^{1/2-\gamma}\left(\theta_n\sqrt{n}\right)^{p-1}=o(n^{p/2})
$$
on the event $A_n\cap\{|z^-+B(n-k)|\> \theta_n\sqrt{n}\}$. Consequently,
\begin{align*}
&\mathbf{E}[|u(z+S(n-k))-u(z^-+B(n-k))|,\tau^{bm}_{z^-}>n-k,A_n,|z^-+B(n-k)|\> \theta_n\sqrt{n}]\\
&\hspace{1cm}=o(n^{p/2}\mathbf{P}(\tau^{bm}_{z^-}>n-k))=o(u(z)).
\end{align*}
Combining this estimate with \eqref{t.h.6} and \eqref{t.h.7} we conclue that \eqref{t.h.4} is 
valid also for $p<1$.
\subsection{Proof of Corollary \ref{c.h}} In view of Theorem \ref{t.h} it suffices to show that $|B(t)|$ under
$\mathbf{P}^{(u)}_0$ is a Markov process with the desired transition kernel.

We first note that \eqref{h-trans.final} implies that
$$
\mathbf{P}^{(u)}_0(B(1)\in dx)/dx=\varkappa u(x)e(1,x)=Cu^2(x)e^{-|x|^2/2}.
$$
Recalling that $u(x)=|x|^pm_1(x/|x|)$, we then obtain
$$
\mathbf{P}^{(u)}_0(|B(1)|\in dr)/dr=Cr^{2p+d-1}e^{-r^2/2}. 
$$
Finally, using the scaling property of $B(t)$, we arrive at the following expression for the entrance law
\begin{equation}
\label{entrance.law}
\mathbf{P}^{(u)}_0(|B(t)|\in dr)/dr=Ct^{-p-d/2}r^{2p+d-1}e^{-r^2/2t}.
\end{equation}

Under $\mathbf{P}^{(u)}_0$, $B(t)$ is a Markov process with the following transition kernel:
\begin{align*}
\frac{\mathbf{P}^{(u)}_0(B(t+h)\in dz|B(t)=x)}{dz}
&=\frac{u(z)}{u(x)}\frac{\mathbf{P}(x+B(h)\in dz,\tau^{bm}_x>h)}{dz}\\
&=\frac{u(z)}{u(x)}b_h(x,z),
\end{align*}
where, according Lemma 1 in \cite{BS97},
$$
b_h(x,z)=\frac{e^{-(|x|^2+|z|^{2})/2h}}{h|x|^{d/2-1}|z|^{d/2-1}}\sum_{j=1}^\infty
I_{a_j}\left(\frac{|x||z|}{h}\right)m_j\left(\frac{x}{|x|}\right) m_j\left(\frac{z}{|z|}\right),
$$
with $a_j=\sqrt{\lambda_j+(d/2-1)^2}$.

Combining this with the entrance law for $B(t)$ we have
\begin{align*}
&\frac{\mathbf{P}^{(u)}_0(B(t+h)\in dz,B(t)\in dx)}{dxdz}=Cu(z)t^{-p-d/2}u(x)e^{-|x|^2/2t}b_h(x,z)
\end{align*}
Integrating over $\{x:|x|=r_1\}$ and $\{z:|z|=r_2\}$ and using orthogonality of eigenfunctions $m_j$,
we obtain
\begin{align*}
&\frac{\mathbf{P}^{(u)}_0(|B(t)|\in dr_1, |B(t+h)|\in dr_2)}{dr_1dr_2}\\
&=Cr_1^{d-1}r_2^{d-1}r_2^pt^{-p-d/2}r_1^pe^{-r_1^2/2t}\frac{e^{-(r_1^2+r_2^{2})/2h}}{hr_1^{d/2-1}r_2^{d/2-1}}
I_{a_1}\left(\frac{r_1r_2}{h}\right).
\end{align*}
Using now the Bayes formula, we obtain
\begin{align*}
&\frac{\mathbf{P}^{(u)}_0\left(|B(t+h)|\in dr_2\big||B(t)|=r_1\right)}{dr_2}\\
&\hspace{1cm}=C \left(\frac{r_2}{r_1}\right)^{p+d/2-1}r_2\frac{e^{-(r_1^2+r_2^{2})/2h}}{h}
I_{a_1}\left(\frac{r_1r_2}{h}\right).
\end{align*}
Noting that $a_1=p+(d/2-1)$ we see that the right hand side is the transition kernel of the $(2p+d)$-dimensional
Bessel process.

By similar calculations one can easily show that the process $|B(t)|$ is markovian under $\mathbf{P}^{(u)}_0$.
Thus, the proof of the corollary is finished.
\section{Convergence of bridges: proof of Theorem \ref{t.bridge}}
In order to simplify our calculations, we shall assume that our random walk is strongly aperiodic.
For every event $B\in\sigma(\{S(k),k\leq nt\})$ we have
\begin{align}
\label{t.bridge.0}
\nonumber
&\mathbf{P}(B|x+S(n)=y,\tau_x>n)\\
\nonumber
&=\frac{\sum_{z\in K}\mathbf{P}(B,x+S(nt)=z,\tau_x>nt)\mathbf{P}(z+S((1-t)n)=y,\tau_z>(1-t)n)}
{\mathbf{P}(x+S(n)=y,\tau_x>n)}\\
&=\mathbf{E}\left[h^{(n)}_{x,y}\left(t,X^{(n)}(t)\right){\bf 1}_B|\tau_x>nt\right],
\end{align}
where
$$
h^{(n)}_{x,y}(t,w)=\frac{\mathbf{P}(\tau_x>nt)\mathbf{P}\left(w\sqrt{n}+S((1-t)n)=y,\tau_{w\sqrt{n}}>(1-t)n\right)}
{\mathbf{P}(x+S(n)=y,\tau_x>n)}.
$$
We are going to show that there exists a bounded, continuous function $h(t,w)$ such that
\begin{equation}
\label{t.bridge.1}
\sup_{w\in K}|h^{(n)}_{x,y}(t,w)-h(t,w)|\to0.
\end{equation}
Let $\widetilde{S}(k)$ denote a random walk, whose increments are independent copies of $-X(1)$. Considering
the path $w\sqrt{n}+S(k),\ k\leq (1-t)\sqrt{n}$ in the reversed time, we have
\begin{align*}
&\mathbf{P}\left(w\sqrt{n}+S((1-t)n)=y,\tau_{w\sqrt{n}}>(1-t)n\right)\\
&\hspace{1cm}=\mathbf{P}\left(y+\widetilde{S}((1-t)n)=w\sqrt{n},\widetilde{\tau}_y>(1-t)n\right).
\end{align*}
Thus, applying Theorem 5 from \cite{DW15}, we conclude that, uniformly in $w\in K$,
\begin{align}
\label{t.bridge.2}
\nonumber
n^{p/2+d/2}\mathbf{P}\left(w\sqrt{n}+S((1-t)n)=y,\tau_{w\sqrt{n}}>(1-t)n\right)&\\
-\varkappa \widetilde{V}(y) H_0 \frac{u(w)}{(1-t)^{p+d/2}}e^{-|w|^2/2(1-t)}&\to0,
\end{align}
where $\widetilde{V}$ is the positive harmonic function for $\widetilde{S}$ and $H_0$
is a norming constant.

Furthermore, combining Theorem 1 and Theorem 6 from \cite{DW15}, we get
$$
\frac{\mathbf{P}(\tau_x>nt)}{\mathbf{P}(x+S(n)=y,\tau_x>n)}\sim
C\frac{t^{-p/2}}{\widetilde{V}(y)}n^{p/2+d/2}.
$$
From this estimate and \eqref{t.bridge.2} we infer that \eqref{t.bridge.1} is valid with
$$
h(t,w)=Ct^{-p/2}(1-t)^{-p/2-d/2}u(w)e^{-|w|^2/2(1-t)}.
$$

Let $g_t:D[0,t]\mapsto\mathbb{R}$ be bounded and continuous. It follows from \eqref{t.bridge.0} and
\eqref{t.bridge.1} that
\begin{equation*}
\mathbf{E}[g_t(X^{(n)})|x+S(n)=y,\tau_x>n]
=(1+o(1))\mathbf{E}[g_t(X^{(n)})h(t,X^{(n)}(t))|\tau_x>nt].
\end{equation*}
Applying now Theorem \ref{t.meander}, we finally get
\begin{equation}
\label{t.bridge.3}
\mathbf{E}[g_t(X^{(n)})|x+S(n)=y,\tau_x>n]\to
\mathbf{E}[g_t(t^{1/2}M_k)h(t,t^{1/2}M_K(1))].
\end{equation}
In other words, we have shown convergence in distribution on $D[0,t]$ for every fixed $t<1$.
Furthermore, the right hand side in \eqref{t.bridge.3} determines all finite dimensional distributions
of the suggested limiting process. Thus, it remains to show tightness on the time interval $[1-\delta,1]$.

Considering the path $x+S(k),\ k\leq n$ in reversed time and using the same arguments, one obtains
for every bounded and continuous $q_t:D[t,1]\mapsto\mathbb{R}$ convergence
\begin{align}
\label{t.bridge.4}
\nonumber
&\mathbf{E}[q_t(X^{(n)})|x+S(n)=y,\tau_x>n]\\
&\hspace{1cm}\to\mathbf{E}[q_t((1-t)^{1/2}M_K)h(1-t,(1-t)^{1/2}M_K(1))].
\end{align}
As a result, we have tightness on $[1-\delta,1]$ for every $1>\delta>0$.
Therefore, the proof of the weak convergence is finished. 

Finally, we note that the continuity of $M_K$ implies also that
$B_K^0$ is continuous as well. Moreover, due to $M_K(0)=0$, we have $B^0_K(0)=B^0_K(1)=0.$\vspace{2mm}\\
\paragraph{\textbf{Acknowledgment}} The first author thanks University of Augsburg, Germany for its hospitality during a visit in July 2015 during which part of this work was completed. 


\begin{thebibliography}{99}
\bibitem{BS97} Banuelos, R. and Smits, R.G.
\newblock Brownian motion in cones.
\newblock \emph{Probab. Theory Related Fields}, {\bf 108}:299-319, 1997.
  
\bibitem{Billing} Billingsley, P.
\newblock {\em Convergence of probability measures.}
\newblock John Wiley and Sons, 1968.

\bibitem{Bolt76} Bolthausen, E.
\newblock On a functional central limit theorem for random walks conditioned to stay positive.
\newblock {\em Ann. Probab.}, {\bf 4}:480-485, 1976.

\bibitem{Bor72} Borovkov, A.A.
\newblock Notes on inequalities for sums of independent random variables.
\newblock \emph{Theory Probab. Appl.}, {\bf 17}:556-557, 1972.

\bibitem{BD06} Bryn-Jones, A. and Doney, R.A.
\newblock A functional limit theorem for random walk conditioned to stay non-negative.
\newblock\emph{J. London Math. Soc. (2)}, {\bf 74}:244-258, 2006.

\bibitem{CC08}Caravenna, F. and Chaumont, L.
\newblock Invariance principles for random walks conditioned to stay positive.
\newblock {\em Ann. Inst. H. Poincare Probab. Statist.,} {\bf 44}:170-190, 2008.

\bibitem{CC13} Caravenna, F. and Chaumont, L.
\newblock An invariance principle for random walk bridges conditioned to stay positive.
\newblock {\em Electron. J. Probab.}, {\bf 18}, Paper no. 60, 2010.

\bibitem{Ch97} Chaumont, L.
\newblock Excursion normalis\'{e}e, m\'{e}andre et pont pour les processus de L\'{e}vy stables.
\newblock {\em Bull. Sci. Math.}, {\bf 121}:377-403, 1997.

\bibitem{DW10} Denisov, D. and Wachtel, V.
\newblock Conditional limit theorems for ordered random walks.
\newblock {\em Electron. J. Probab.}, {\bf 15}:292--322, 2010.

\bibitem{DW15} Denisov, D. and Wachtel, V.
\newblock Random walks in cones.
\newblock {\em Ann. Probab.}, {\bf 43}:992-1044, 2015.

\bibitem{Don85} Doney R.A.
\newblock Conditional limit theorems for asymptotically stable random walks.
\newblock{\em Probab. Theory Relat. Fields,} {\bf 70}:351-360, 1985.

\bibitem{Feierl12} Feierl, T.
\newblock The height of watermelons with wall. 
\newblock {\em J. Phys. A}, {\bf 45}, Paper no. 9, 2012. 

\bibitem{Feierl13} Feierl, T.
\newblock The height and range of watermelons without wall.
\newblock {\em  European J. Combin.}, {\bf 34}:138-154, 2013.

\bibitem{Gar09} Garbit, R.
\newblock Brownian motion conditioned to stay in a cone.  
\newblock {\em J. Math. Kyoto Univ.}, {\bf 49}:573-592, 2009.

\bibitem{Gar11} Garbit, R.
\newblock A central limit theorem for two-dimensional random walks in a cone.
\newblock {\em Bulletin de la SMF}, {\bf 139}:271-286, 2011.

\bibitem{Grab99} Grabiner, D.J.
\newblock Brownian motion in a Weyl chamber, non-colliding particles, and random matrices.
\newblock{\em Ann. Inst. H. Poincare Probab. Statist.}, {\bf 35}(2):177-204, 1999.

\bibitem{Shim91} Shimura, M.
\newblock A limit theorem for two-dimensional random walk conditioned to stay in a cone.
\newblock {\em Yokohama Math. J.}, {\bf 39}:21-36, 1991.

\bibitem{Sp76} Spitzer, F.
\newblock {\em Principles of random walk}, 2nd edition.
\newblock Springer, New York, 1976.

\bibitem{Stone65} Stone, C.
\newblock On local and ratio limit theorems. 
\newblock Proceedings of the Fifth Berkeley Symposium on Mathematical Statistics and Probability, 
Volume 2: Contributions to Probability Theory, Part 2, 217--224, University of California Press, Berkeley, Calif., 1967.
\end{thebibliography}
\end{document}